%% file: qcerf2tor.tex
\documentclass[10pt]{amsart}  
\usepackage{graphicx}
\usepackage[pdftex, plainpages=false]{hyperref}
\usepackage{color}
\newtheorem{theorem}[subsection]{Theorem}
\newtheorem{proposition}[subsection]{Proposition}
\newtheorem{lemma}[subsection]{Lemma}
\newtheorem{corollary}[subsection]{Corollary}
\newtheorem{claim}[subsection]{Claim}
\theoremstyle{definition}
\newtheorem{definition}[subsection]{Definition}
\newtheorem{remark}[subsection]{Remark}
\newtheorem{case}{Case}

\newcommand{\nbd}[2]{\mathcal{N}_{#2}({#1})}  
\newcommand{\mc}[1]{\mathcal{#1}}
\newcommand{\HH}{\mathbb{H}}
\newcommand{\ZZ}{\mathbb{Z}}
\newcommand{\co}{\colon\thinspace}

\newcommand{\Label}{\mathrm{Label}}

\begin{document}
\title[Separation of Relatively Quasiconvex Subgroups]{Separation of Relatively
  Quasiconvex Subgroups}


\author[J. F. Manning]{Jason Fox Manning}
\address{Department of Mathematics, SUNY at Buffalo,
Buffalo, NY 14260-2900}
\email{j399m@buffalo.edu}

\author[E. Mart\'inez-Pedroza]{Eduardo Mart\'inez-Pedroza}
\address{Department of Mathematics \& Statistics, McMaster University, Hamilton, ON, L8S 4K1, Canada.}
\email{emartinez@math.mcmaster.ca}

\begin{abstract}
Suppose that all hyperbolic groups are residually finite.  The
following statements follow:  In relatively hyperbolic groups with
peripheral structures consisting of finitely generated nilpotent
subgroups, relatively quasiconvex subgroups are separable; Geometrically finite
subgroups of non-uniform lattices in rank one symmetric spaces are
separable; Kleinian groups are subgroup separable.  We also show that
LERF for finite volume hyperbolic $3$--manifolds would follow from
LERF for closed hyperbolic $3$--manifolds.

The method is to reduce, via combination and filling theorems, the
separability of a relatively quasiconvex subgroup of a relatively hyperbolic
group $G$ to the separability of a quasiconvex subgroup of a
hyperbolic quotient $\bar{G}$.  A result of Agol, Groves, and Manning
is then applied.
\end{abstract}
\maketitle

\section{Main Results}

A subgroup $H$ of a group $G$ is called \emph{separable} if for any $g
\in G \setminus H$ there is a homomorphism $\pi$ onto a finite group
such that $\pi (g) \not \in \pi (H)$.  A group is called
\emph{residually finite} if the trivial subgroup is separable, and a
group is called \emph{subgroup separable} or \emph{LERF} if every
finitely generated subgroup is separable.  For example, Hall showed
that free groups are
LERF in ~\cite{Hal49}.
It follows from a theorem of Mal$'$cev  \cite{Mal83} that polycyclic
(and in particular finitely generated nilpotent) groups are LERF.  A group is called \emph{slender} if every subgroup is
finitely generated.  Polycyclic groups are also slender, by a result
of Hirsch \cite{Hi37}.

Given a relatively hyperbolic group with peripheral structure
consisting of LERF and slender subgroups, we study separability of
relatively quasiconvex subgroups.  This is connected, via filling
constructions, to residual finiteness of
hyperbolic groups.  It is not known whether all word-hyperbolic groups are residually
finite.  Consequences of a positive or negative answer to this question
have been explored by several authors; see for example 
~\cite{AGM08,  KW00, LN91, Mi06, Os06-1, Wi02}.  In particular, the main result of 
~\cite{AGM08} is the following.
\begin{theorem} \cite{AGM08} \label{thm:AGM}
If all hyperbolic groups are residually finite, then every quasiconvex subgroup of a hyperbolic group is separable.
\end{theorem}
We extend this result, answering a question in \cite{AGM08}, as follows:
\begin{theorem}\label{t:c1}
Suppose that all hyperbolic groups are residually finite.  If $G$ is a
torsion free relatively hyperbolic group with peripheral structure
consisting of subgroups which are LERF and slender, then any
relatively quasiconvex subgroup of $G$ is separable.
\end{theorem}

This extension has some interesting corollaries.
A \emph{pinched Hadamard manifold} is a simply connected Riemannian
manifold with pinched negative curvature.
In \cite{Bo95}, Bowditch gave several equivalent definitions of
geometrical finiteness for discrete subgroups of the isometry group of
a pinched Hadamard manifold, generalizing the
notion of geometrical finiteness in Kleinian groups.  
The next theorem summarizes some useful facts about these groups.
(Statement \eqref{hk:rh} can be found in \cite{Fa98} or \cite{BO99};
statement \eqref{margulis} follows from the Margulis lemma; and
statement \eqref{hk:qc} is \cite[Corollary 1.3]{Hr08}.)

\begin{theorem} \cite{Fa98, Hr08, BO99}\label{thm:HK}
Let $X$ be a pinched Hadamard manifold and let $G$ be a geometrically finite subgroup of $\text{Isom}(X)$.
\begin{enumerate}
\item\label{hk:rh} $G$ is relatively hyperbolic, relative to a
  collection of representatives of conjugacy classes of maximal
  parabolic subgroups. 
\item\label{margulis} Maximal parabolic subgroups of $G$ are virtually
  nilpotent. 
\item\label{hk:qc} A subgroup $H$ of $G$ is relatively
  quasiconvex if and only if $H$ is geometrically finite. 
\end{enumerate}
\end{theorem}
Rank one symmetric spaces are pinched Hadamard manifolds.  We
therefore have the following corollary of Theorem \ref{t:c1}: 
\begin{corollary}
Suppose that all hyperbolic groups are residually finite.  Let
  $G$ be a discrete, geometrically finite subgroup of the
  isometry group of a rank one symmetric space. (For example, $G$
  could be a lattice.)  All the geometrically finite subgroups of $G$ are separable.
\end{corollary}

In case the symmetric space is $\HH^3$, more can be said (see Section
\ref{sec:3mfd} for the proof).
\begin{corollary}\label{cor:Kleinian}
  If all hyperbolic groups are residually finite,
  then all finitely generated Kleinian groups are LERF\@.
\end{corollary}

Briefly, Theorem \ref{t:c1} is proved by combining one of
Mart\'inez-Pedroza's combination theorems in \cite{MP07} with
Theorem \ref{thm:AGM} and the Dehn filling technique of \cite{GM06,Os06-1}.
We next give a more detailed discussion.  


\begin{definition}   
A relatively quasiconvex subgroup $H$ of $G$ is called \emph{fully  quasiconvex} if for any subgroup $P \in \mathcal{P}$ and any $f\in
G$, either $H\cap P^f$ is finite or $H\cap P^f$ is a finite index
subgroup of $P^f$.  (Here $P^f = fPf^{-1}$.)
\end{definition}
Using the work in \cite{MP07}, we show the following.
\begin{theorem} \label{thm:fully-quasiconvex}
Let $G$ be a group hyperbolic relative to a collection of slender and LERF subgroups.
Suppose that $Q$ is a relatively quasiconvex subgroup of $G$ and $g$ is an element of $G$ not in $Q$.  
Then there exists a fully quasiconvex subgroup $H$ which contains $Q$ and does not contain $g$.
\end{theorem}

\begin{remark}\label{r:priorart}
  In case $G$ is a finite volume   hyperbolic $3$--manifold group and $Q$ is the fundamental group of a   quasi-fuchsian surface, Theorem   \ref{thm:fully-quasiconvex} can be proved using geometric   arguments like those in \cite{CL99,CL01}.   Such geometric arguments were applied (in a   different way) to   separability questions in \cite{ALR} (cf. \cite{Wi06}).
\end{remark}

Using the work in ~\cite{AGM08} and Theorem~\ref{thm:fully-quasiconvex}, we prove the following separation result of relatively quasiconvex subgroups by maps onto word hyperbolic groups.  (A part of this theorem can be interpreted as saying that $G$ is ``quasiconvex extended residually hyperbolic''.)

\begin{theorem} \label{thm:1}
Let $G$ be a torsion free group hyperbolic relative to a collection of slender and LERF subgroups.
For any relatively quasiconvex subgroup $Q$ of $G$ and any element $g \in G$ such that $g \not \in Q$, 
there is a fully quasiconvex subgroup $H$ of $G$, and a surjective
homomorphism $\pi \co G \longrightarrow \bar{G}$  such that 
\begin{enumerate}
\item\label{1qinh} $Q<H$,
\item\label{1wh} $\bar{G}$ is a word-hyperbolic group, 
\item\label{1qc} $\pi (H)$ is a quasiconvex subgroup of $\bar{G}$, 
\item\label{1sep} $\pi (g) \not \in \pi (H)$.
\end{enumerate}
\end{theorem}
We can prove Theorem \ref{t:c1} from Theorem \ref{thm:1} as follows:
\begin{proof}[Proof of Theorem \ref{t:c1}]
  Let $Q<G$ be relatively quasiconvex, and let $g\in G\setminus Q$.  By
  Theorem \ref{thm:1} there is a fully quasiconvex $H<G$ containing
  $Q$ but not $g$, and a quotient $\pi\co G\to K$ so that $\pi(g)\notin
  \pi(H)$, $K$ is hyperbolic, and $\pi(K)$ is quasiconvex.  

  Assuming all hyperbolic groups are residually finite, Theorem
  \ref{thm:AGM} implies that there is a finite group $F$ and a
  quotient $\phi\co K\to F$ so that $\phi(\pi(g))\notin
  \phi(\pi(H))$.  Since $\phi(\pi(H))$ contains $\phi(\pi(Q))$, the
  map $\phi\circ\pi$ serves to separate $g$ from $Q$.
\end{proof}
\begin{remark}
  The torsion-free hypotheses
  in Theorems \ref{thm:1} and \ref{t:c1} 
  are not really necessary.  We sketch the necessary changes to our
  argument in Appendix \ref{app:torsion}.
  If one is primarily interested in Theorem
  \ref{t:c1} in the special case of virtually polycyclic peripheral
  subgroups, we have the following simple argument, pointed out to us by the
  referee:

  Let $G$ be a relatively hyperbolic group, relative to a collection
  $\mathcal{P}=\{P_1,\ldots,P_m\}$ of virtually polycyclic subgroups.  An
  easy argument shows that each $P_i$ contains a finite index normal 
  subgroup which is torsion free.  Moreover, $G$ contains only finitely many finite
  order non-parabolic elements, up to conjugacy \cite[Theorem 4.2]{Os06}. 
  It then follows from Osin's version of the relatively hyperbolic
  Dehn filling theorem that there is a filling 
  $G \stackrel{\pi}{\to} G(N_1',\ldots, N_m')$, so that $G(N_1',\ldots, N_m')$ is
  hyperbolic, and no non-trivial torsion element of $G$ is in the kernel of
  $\pi$.  Assuming hyperbolic groups are residually finite, the group
  $G(N_1',\ldots,N_m')$ has a torsion-free subgroup $S$ of finite
  index.  The preimage $G_0=\pi^{-1}(S)$ is a torsion-free finite index
  subgroup of $G$.  Again under the assumption that hyperbolic groups
  are residually finite, Theorem \ref{t:c1} implies that $G_0$ is
  QCERF, and so (applying Corollary 2.2 below) $G$ must also be
  QCERF\@.
\end{remark}

The paper is organized as follows.  In Section \ref{sec:def} we give a
definition of relatively hyperbolic group which suits the purposes of
this paper and is equivalent to the more standard definitions in the
literature.  In Section \ref{sec:comb} we recall a combination theorem
for relatively quasiconvex subgroups from \cite{MP07} and prove
Theorem~\ref{thm:fully-quasiconvex}.  In Section \ref{sec:fill} we
recall some definitions and results on fillings of relatively
hyperbolic groups and prove Theorem~\ref{thm:1}.  In Section
\ref{sec:3mfd}, we give two applications to separability questions on
hyperbolic $3$--manifolds: Corollary~\ref{cor:Kleinian} and
Proposition \ref{c:fv}.  In Appendix \ref{app:qcequiv}, we prove a
result we need on the equivalence of various definitions of relative
quasiconvexity, and in Appendix \ref{app:torsion}, we sketch how to
prove our results in the presence of torsion.

\section{Preliminaries}
\subsection{Separability}
Let $G$ be a group.
Recall that the \emph{profinite topology} on $G$ is the smallest
topology on $G$ in which all finite index subgroups and their cosets
are closed.  The group $G$ is residually finite if and only if this
topology is Hausdorff.  A subgroup $H$ is separable if and only if it
is a closed subset of $G$, with this topology.

Given a subgroup $G_0<G$, one can ask whether the profinite topology
on $G_0$ coincides with the topology induced by the profinite topology
on $G$.  In general, the topologies are quite different, but in case
$G_0$ is finite index, the topologies coincide.
In particular, we have:
\begin{lemma}
  Let $G_0<G$ be a finite index subgroup, and let $C_0\subseteq G_0$.
  The following conditions are equivalent:
  \begin{enumerate}
  \item\label{proclosed} $C_0$ is closed in the profinite topology on $G_0$.
  \item\label{bigclosed} $C_0$ is closed in the profinite topology on $G$.
  \item\label{subclosed} $C_0 = C\cap G_0$ for some set $C$ which is closed in the
    profinite topology on $G$.
  \end{enumerate}
\end{lemma}
\begin{proof}
  Suppose $C_0$ is closed in $G_0$.  Since the finite index subgroups of
  $G_0$ and their cosets generate the topology on $G_0$, we can write
  $C_0$   as a finite union of arbitrary intersections of cosets
\begin{equation}\label{unint1}
  C_0 = \bigcup_{i=1}^n \bigcap_{i\in I_j} g_i K_i 
\end{equation}
  where every $K_i$ is a finite index subgroup of $G_0$.  But since
  $G_0$ is finite index in $G$, each of the $K_i$ appearing in equation
  \eqref{unint1} is also finite index in $G$.  This $C_0$ is a finite
  union of intersections of closed sets in the profinite topology on
  $G$, so $C_0$ is closed in $G$.  Thus condition \eqref{proclosed}
  implies condition \eqref{bigclosed}.

  Trivially, condition \eqref{bigclosed} implies condition
  \eqref{subclosed}.  To see that condition \eqref{subclosed} implies
  condition \eqref{proclosed}, we first establish the following.
  \begin{claim}
    If $g\in G$, and $K<G$ is finite index, then $gK\cap G_0$ is
    closed in the profinite topology on $G_0$.
  \end{claim}
  \begin{proof}
    Let $K_0 = K\cap G_0$.  
    Since $K_0$ is finite index in  $K$, the subgroup $K$ is a finite union of
    cosets $K = \cup_{i=1}^p g_i K_0$, and so $gK = \cup_{i=1}^p gg_i
    K_0$ is as well.  
    Since a coset of $K_0$ in $G$ either lies inside $G_0$ or in its complement,
    it follows that $gK\cap G_0$ is a finite union of cosets of $K_0$ in $G_0$.
    Since $K_0$ is finite index in $G_0$, these cosets are closed in the profinite topology on $G_0$ and therefore $gK\cap G_0$ is closed as well.
  \end{proof}
  Suppose then that $C_0 = C\cap G_0$ for $C$ closed in $G$.  We have
    $C = \bigcup_{i=1}^n \bigcap_{i\in I_j} g_i K_i$
  where now the $K_i$ are arbitrary finite index subgroups of $G$, and
  the $g_i$ are arbitrary elements of $G$.  We thus have
  \begin{equation}\label{unint}
    C_0 = C\cap G_0 = \bigcup_{i=1}^n \bigcap_{i\in I_j} (g_i K_i)\cap G_0.
  \end{equation}
  By the claim, every $(g_i K_i)\cap G_0$ appearing in equation
  \eqref{unint} is closed in $G_0$, and so
  $C_0$ is also closed in the profinite topology on $G_0$.
\end{proof}
\begin{corollary}\label{l:finindex}
  Let $H<G$ be a pair of groups, and let $G_0$ be
  a finite index subgroup of $G$.  Let $H_0 = H\cap G_0$.  The
  subgroup $H$ is separable in $G$ if and only if $H_0$ is separable
  in $G_0$.
\end{corollary}
Another immediate corollary is that LERFness and QCERFness are
commensurability invariants.

\subsection{Relative hyperbolicity}\label{sec:def}

The notion of relative hyperbolicity has been studied by several authors with different equivalent definitions. The definition in this subsection is based on the work by D.\ Osin in~\cite{Os06}.  Let $G$ be a group, $\mathcal{P}$ denote a collection  of subgroups $\{ P_1, \dots , P_m \}$, and $S$ be a finite generating set which is assumed to be symmetric, i.e, $S=S^{-1}$.  Denote by $\Gamma(G, \mathcal{P}, S)$ the Cayley graph of $G$ with respect to the generating set $S\cup \bigcup \mathcal{P}$.  If $p$ is a path between vertices in $\Gamma(G, \mathcal{P}, S)$, we will refer to its initial vertex as $p_-$, and its terminal vertex as $p_+$.  The path $p$ determines a word $\Label(p)$ in the alphabet $S\cup \bigcup \mathcal{P}$ which represents an element $g$ so that $p_+ = p_- g$.

\begin{definition}[Weak Relative Hyperbolicity]
The pair $(G, \mathcal{P})$ is \emph{weakly relatively hyperbolic} if
there is an integer $\delta \geq 0$ such that 
$\Gamma(G, \mathcal{P},S)$ is $\delta$--hyperbolic.  We may also say
that $G$ is \emph{weakly relatively hyperbolic, relative to}
$\mathcal{P}$.  
\end{definition}

\begin{definition}[~\cite{Os06}]
Let $q$ be a combinatorial path in the Cayley graph $\Gamma(G,
\mathcal{P}, S)$.  Sub-paths of $q$ with at least one edge are called
\emph{non-trivial}.  For $P_i \in \mathcal{P}$, a \emph{ $P_i$--component} of
$q$ is a maximal non-trivial sub-path $s$ of $q$ with $Label(s)$ a word
in the alphabet $P_i$.  When we don't need to specify the index $i$,
we will refer to $P_i$--components as
\emph{$\mathcal{P}$--components}.

Two $\mathcal{P}$--components $s_1$, $s_2$  are \emph{connected}
if the vertices of $s_1$ and $s_2$ belong to the same left coset of
$P_i$ for some $i$.  A $\mathcal{P}$--component $s$ of $q$ is 
\emph{ isolated} if it is not connected to a different
$\mathcal{P}$--component of $q$.  The path $q$ is 
\emph{ without backtracking} if every $\mathcal{P}$--component of $q$ is isolated.

A vertex $v$ of $q$ is called \emph{ phase} if it is not an interior
vertex of a $\mathcal{P}$--component $s$ of $q$.  Let $p$ and $q$ be
paths between vertices in $\Gamma(G, \mathcal{P}, S)$.  The paths $p$
and $q$ are \emph{ $k$--similar} if
\[ \max\{dist_{S}(p_-,q_-), dist_{S}(p_+, q_+)  \} \leq k,\]
where $dist_S$ is the metric induced by the finite generating set $S$
(as opposed to the metric in $\Gamma(G, \mathcal{P}, S)$).
\end{definition}
\begin{remark}
A geodesic path $q$ in $\Gamma(G, \mathcal{P}, S)$ is without
backtracking, all $\mathcal{P}$--components of $q$ consist of a single
edge, and all vertices of $q$ are phase. 
\end{remark}

\begin{definition}[Bounded Coset Penetration (BCP)] \label{BCP}
The pair $(G, \mathcal{P})$ satisfies the \emph{BCP property} if for
any $\lambda \geq 1$, $c \geq 0$, $k\geq 0$, there exists an integer
$\epsilon(\lambda, c, k) > 0$ such that for $p$ and $q$ any two $k$--similar
$(\lambda, c)$--quasi-geodesics in $\Gamma(G, \mathcal{P}, S)$ without
backtracking, the following conditions hold:
\begin{enumerate}
\item[(i.)] The sets of phase vertices of $p$ and $q$ are contained in the closed
$\epsilon(\lambda, c, k)$--neighborhoods  of each other, with respect to the metric $dist_S$.
\item[(ii.)] If $s$ is any $\mathcal{P}$--component of $p$ such that $dist_S(s_-, s_+) > \epsilon(\lambda, c, k)$, then
there exists a $\mathcal{P}$--component $t$ of $q$ which is connected to $s$.
\item[(iii.)] If $s$ and $t$ are connected $\mathcal{P}$--components of $p$ and $q$ respectively,
then \[\max\{ dist_S(s_-, t_-), dist_S(s_+, t_+) \} \leq  \epsilon(\lambda, c, k).\]
\end{enumerate}
\end{definition}
\begin{remark}
Our definition of the BCP property corresponds to the conclusion of Theorem 3.23 in \cite{Os06}. 
\end{remark}
\begin{definition}[Relative Hyperbolicity] \label{defn:rel_hyp}
The pair $(G, \mathcal{P})$ is \emph{relatively hyperbolic} if the
group $G$ is weakly relatively hyperbolic relative to $\mathcal{P}$ and the pair
$(G,\mathcal{P})$ satisfies the Bounded Coset Penetration property.  If
$(G, \mathcal{P})$ is relatively hyperbolic then
we say $G$ is \emph{relatively hyperbolic, relative to}
$\mathcal{P}$; if there is no ambiguity, we just say that
$G$ is relatively hyperbolic. 
\end{definition}
\begin{remark}
Definition \ref{defn:rel_hyp} given here is equivalent to Osin's
\cite[Definition 2.35]{Os06} for finitely generated groups: To see
that Osin's definition implies \ref{defn:rel_hyp}, apply
\cite[Theorems 3.23]{Os06}; to see that \ref{defn:rel_hyp} implies
Osin's definition, apply \cite[Lemma 7.9 and Theorem 7.10]{Os06}.  For
the equivalence of Osin's definition and the various other definitions
of relative hyperbolicity see \cite{Hr08} and the references therein.

The definition of relative hyperbolicity is independent of finite
generating set $S$.
\end{remark}

\section{Combination of Parabolic and quasiconvex Subgroups}\label{sec:comb}

In this section, $G$ will be relatively hyperbolic, relative to a finite
collection of subgroups $\mathcal{P}$, and $S$ will be a finite
generating set for $G$.  Denote by $\Gamma(G,
\mathcal{P}, S)$ the Cayley graph of $G$ with respect to the
generating set $S\cup \bigcup \mathcal{P}$.

\subsection{Parabolic and Quasiconvex Subgroups}
\begin{definition}
The \emph{peripheral} subgroups of $G$ are the elements of $\mathcal{P}$.
A subgroup of $G$ is called \emph{parabolic} if it can be conjugated
into a peripheral subgroup.
\end{definition}

\begin{proposition}\cite[Proposition 2.36]{Os06} \label{prop:parabolic-subgroups}
The following conditions hold.
\begin{enumerate}
\item For any $g_1$, $g_2\in G$, the intersection $P_i^{g_1} \cap P_j^{g_2}$ is finite unless $i = j$.
\item The intersection $P_i^{g} \cap P_i$ is finite for any $g \not \in P_i$.
\end{enumerate}
In particular, if $Q$ is a subgroup of $G$, then any infinite maximal parabolic subgroup of $Q$ is of the form $Q\cap P_i^f$ for some $f \in Q$ and $P_i \in \mathcal{P}$.
\end{proposition}

\begin{definition}\cite[Definition 4.9]{Os06} \label{def:O}
A subgroup $Q$ of $G$ is called \emph{quasiconvex relative to
  $\mathcal{P}$} (or simply \emph{relatively quasiconvex} when the
collection $\mathcal{P}$ is fixed) if there exists a constant $\sigma
\geq 0$ such that the following holds: Let $f$, $g$ be two elements of
$Q$, and $p$ an arbitrary geodesic path from $f$ to $g$ in the Cayley
graph $\Gamma(G, \mathcal{P}, S)$.  For any vertex $v
\in p$, there exists a vertex $w \in Q$ such that $ dist_S(v,w) \leq
\sigma,$ where $dist_S$ is the word metric induced by $S$.
\end{definition}

\begin{remark}
For more on different definitions of relative quasiconvexity in the literature, see
Appendix \ref{app:qcequiv}.
\end{remark}

\begin{theorem}\cite[Theorem 9.1]{Hr08} \label{prop:ParabolicClasses}
Let $Q$ be a finitely generated relatively quasiconvex subgroup of
$G$.  The number of infinite maximal parabolic subgroups of $Q$ up
to conjugacy in $Q$ is finite.  Furthermore, if $\mathcal{O}$ is a set
of representatives of these conjugacy classes, then $Q$ is relatively
hyperbolic, relative to $\mathcal{O}$. 
\end{theorem}
\begin{remark}
In \cite{Hr08}, an extended definition of relative hyperbolicity is used
which includes some countable but non-finitely generated groups. 
Using this extended definition, 
the assumption of finite generation
in Theorem~\ref{prop:ParabolicClasses} is superfluous.  

We note that in case all the peripheral subgroups are slender,
relatively quasiconvex subgroups are necessarily finitely generated (see
\cite[Corollary 9.2]{Hr08}).
\end{remark}

\subsection{Combination of Quasiconvex Subgroups}
For $g \in G$, $|g|_S$ denotes the distance from $g$ to the identity element in the word metric induced by $S$. 

\begin{theorem} \cite[Theorem 1.1]{MP07} \label{thm:2}
Let $Q$ be a relatively quasiconvex subgroup of $G$, and let $P$ be a
maximal parabolic subgroup of $G$.  Suppose that $P^f=P_i$ for some
$P_i \in \mathcal{P}$ and $f \in G$.

There are constants $C=C(Q, P)\geq 0$ and $c=c(Q,P)\geq 0$ with the
following property.  Suppose $D \geq C$ and $R$ is a subgroup of $P$
such that
\begin{itemize}
\item $P\cap Q < R$, and
\item $|g|_S >D$ for any element $g \in R\setminus Q$.
\end{itemize}
It follows that:
\begin{enumerate}
\item\label{i:qc}  The subgroup $H=\langle Q\cup R \rangle$ is relatively quasiconvex and the natural map  $Q\ast_{Q\cap R} R \longrightarrow H$ is an isomorphism.

\item\label{i:pc} Every parabolic subgroup of $H$ is conjugate in $H$ to a parabolic subgroup of $Q$ or $R$.

\item\label{i:ge} For any $g \in H$, either $g \in Q$, or any geodesic from $1$ to $g$ in the relative Cayley graph $\Gamma(G, \mathcal{P}, S)$ has 
at least one $P_i$-component $t$ such that $|t|_S > D-c$.
\end{enumerate}
\end{theorem}
\begin{proof}
Conclusions \eqref{i:qc} and \eqref{i:pc} rephrase ~\cite[Theorem 1.1]{MP07}.  The proof of conclusion~\eqref{i:ge} is divided into two cases:  $P \in \mathcal{P}$ and $P \notin \mathcal{P}$.

\begin{case}\label{inP} $P\in \mathcal{P}$.
\end{case}

We summarize part of the argument for~\cite[Theorem 1.1]{MP07} for
conclusions~\eqref{i:qc} and~\eqref{i:pc}; we then explain how conclusion
\eqref{i:ge} follows in this case.

Let $g\in Q\ast_{Q\cap R} R\setminus Q$.  The element $g$ has a
\emph{normal form}
\begin{equation}\label{eq:normalform} g = g_1h_1 \dots
  g_kh_k  \end{equation}
where $g_j \in Q\setminus Q\cap R$ for $1< j \leq k$,  $h_j \in R
\setminus Q\cap R$ for $1\leq j < k$, either $g_1=1$ or $g_1 \in
Q\setminus Q\cap R$, and either $h_k=1$ or $h_k \in R \setminus Q\cap
R$.  We use the normal form 
to produce a path $o$ in $\Gamma(G, \mathcal{P}, S)$ from $1$ to the
image of $g$ by the natural map  $Q\ast_{Q\cap R} R \longrightarrow
H$ as follows.  For each $j$ between $1$ and $k$, let $u_j$ be a
geodesic path in $\Gamma(G, \mathcal{P},S)$ from $g_1h_1\cdots h_j$ to
$g_1h_1\cdots h_jg_j$ (so that $\Label(u_j)$ represents $g_j$).
Similarly, let $v_j$ be a geodesic path from $g_1h_1\cdots g_{j-1}$ to
$g_1h_1\cdots g_{j-1}h_j$ (so that $\Label(v_j)$ represents $h_j$).
A path $o$ from $1$ to $g$ 
in $\Gamma(G, \mathcal{P}, S)$ is given by
\begin{eqnarray*} o = u_1  v_1  \dots  u_k  v_k.   \end{eqnarray*}
(See Figure~\ref{fig:02}.) 

\begin{figure}[ht]
\begin{center}
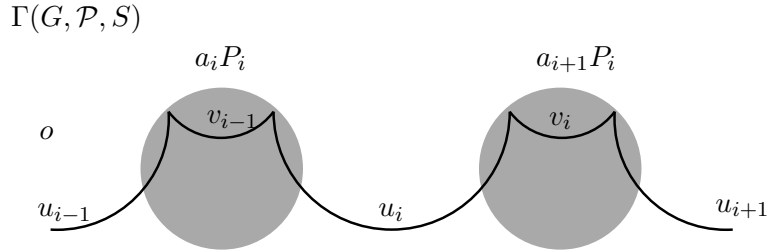
\end{center}
\caption{Part of the polygonal path $o$ in $\Gamma(G, \mathcal{P},
  S)$.  For each $i$, the $P_i$--component of $o$ containing the
  subsegment $v_j$ is long with respect to the $S$--metric.  This implies
  that the path $o$ is a quasi-geodesic with different end-points. } 
\label{fig:02}
\end{figure}

Each subsegment $v_j$ is part of a $\mathcal{P}$-component $t_j$ of
the path $o$.   Let $D$ be as in the hypothesis of the theorem.
The penultimate inequality in the proof of Claim 2 in
the proof of Lemma 5.1 of \cite{MP07} is
\begin{equation}  |t_j|_S = dist_S( (t_j)_-, (t_j)_+ ) > D - 2M(P, Q, \sigma ) ,  \end{equation} 
where $\sigma$ is the quasiconvexity constant for $Q$, 
and $M(P, Q, \sigma)$ is the constant provided by
\cite[Lemma 4.2]{MP07} for the subgroups $Q$, $P$, and the constant
$\sigma$.
Let $\eta$ be the constant from Proposition 3.1 of \cite{MP07}.  If 
$D - 2M(P, Q, \sigma ) > \eta$,  
then $o$ is a $(\lambda, 0)$--quasi-geodesic with distinct
endpoints.  Let $C=\eta + 2M(P, Q, \sigma )$. 

It follows from the argument just sketched that if $D > C$, then the
natural map 
$Q\ast_{Q\cap R} R \longrightarrow H$ is an isomorphism.  
It can further be shown that $H$ is a
relatively quasiconvex subgroup and that the parabolic subgroups of $H$ are
conjugate into $Q$ or $R$ by elements of $H$ (See \cite[Lemmas 5.2 and
5.3]{MP07}for details.).  In other words, parts \eqref{i:qc} and
\eqref{i:pc} hold for $P\in \mathcal{P}$ and $C$ as above.

If $g \in  H \setminus Q$ and $p$ is a geodesic from $1$ to $g$ in
$\Gamma(G, \mathcal{P}, S)$, then the $(\lambda, 0)$--quasi-geodesic
$o$ and the geodesic $p$ are $0$--similar.  Since $o$ has a
$\mathcal{P}$--component of $S$--length at least $D - 2M(P, Q, \sigma
)$, the Bounded Coset Penetration property  (Definition~\ref{BCP})
implies that $p$ has a $\mathcal{P}$--component of $S$--length at
least $D- 2M(P, Q, \sigma )-2\epsilon(0, \lambda, 0) $. 
We have verified \eqref{i:ge} of the Theorem for $c$ equal
to
\[ c(Q,P) = 2M(P,Q,\sigma) +2\epsilon(0,\lambda,0) \]
in the special case that $f=1$ and $P\in \mathcal{P}$.

\begin{case} $P\notin \mathcal{P}$, but $P^f\in \mathcal{P}$ for some $f\in
  G\setminus P$.
\end{case}
Since $P^f\in \mathcal{P}$ and $Q^f$ is relatively quasiconvex, by~\cite[Theorem
  1.1]{MP07} and Case~\ref{inP}, all three conclusions of
Theorem~\ref{thm:2} hold for $Q^f$ and $P^f$ and some constants $C' =
C(Q^f,P^f)>0$ and $c' = c(Q^f,P^f)>0$.   Define   
\[ C=C' + 2|f|_S + 3\epsilon (1, 0, |f|_S) ,\] and \[ c = c' +2|f|_S + 2\epsilon (1, 0, |f|_S),\] 
where $\epsilon (1, 0, |f|_S)$ is the constant of Definition~\ref{BCP}
on the Bounded Coset Penetration property.  Now we show that 
the theorem holds for the subgroups $P$ and $Q$, and the constants $C$ and
$c$.  Let $R$ be a subgroup of $Q$ satisfying the hypothesis of the
theorem for a constant $D>C$. 

If $r \in R\setminus Q$, then $|r|_S\geq D$, by hypothesis.  It 
follows that
\[  |r^f|_S \geq D - 2|f|_S \geq C'. \]
We therefore have:
\begin{enumerate}
\item The subgroup $H^f=\langle Q^f\cup R^f \rangle$ is relatively
  quasiconvex and the natural map  $Q^f \ast_{(Q\cap R)^f} R^f
  \longrightarrow H^f$ is an isomorphism.  Since relative quasiconvexity is
  preserved by conjugation, $H=\langle Q\cup R \rangle$ is relatively
  quasiconvex.  Obviously the map  $Q\ast_{Q\cap R} R \longrightarrow
  H$ is also an isomorphism.  In other words, conclusion \eqref{i:qc}
  holds for $Q$ and $P$ and the constant $C$.

\item Every parabolic subgroup of $H^f$ is conjugate in $H^f$ to a
  parabolic subgroup of $Q^f$ or $R^f$.  Parabolicity is preserved
  under conjugation, so the same property (conclusion \eqref{i:pc} of
  the theorem) holds for the subgroups $Q$, $R$, and $H$, and the
  constant $C$.

\item For any $h \in H^f$, either $h \in Q^f$, or any geodesic from
  $1$ to $h$ in the relative Cayley graph $\Gamma(G, \mathcal{P}, S)$
  has at least one $P_i$--component $t$ such that $|t|_S >
  D-2|f|_S-c'$.
\end{enumerate}

It remains to see why conclusion \eqref{i:ge} of the Theorem holds with the chosen constant $c$.
Let $g \in H \setminus Q$ and let $p$ be a geodesic from $1$ to $g$ in
$\Gamma(G, \mathcal{P}, S)$.  We must show that $p$ has a
$P_i$--component of $S$--length at least $D-c$.  Let $q$ be a geodesic
from $1$ to $fgf^{-1}$.  Since $fgf^{-1}$ belongs to $H^f \setminus
Q^f$, the geodesic $q$ has a $P_i$--component $u$ of $S$--length 
at least $D-2|f|_S -  c'$.
Let $r$ be the geodesic starting at $f$, with the same label as
$p$.  Thus $r$ joins $f$ to $fg$, and
the geodesics $q$ and $r$ are $|f|_S$--similar (see Figure~\ref{fig:01}).
\begin{figure}[ht]
\begin{center}
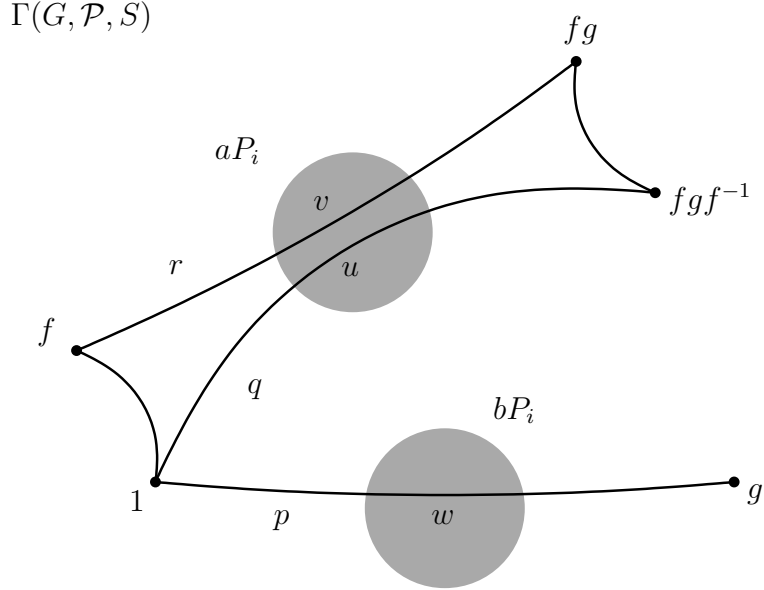
\end{center}
\caption{The geodesics $q$ and $r$ are $|f|_S$--similar and $q$ contains a large $\mathcal{P}$--component $u$.  By the BCP-property, $r$ has a large $\mathcal{P}$--component $v$.  Since $p$ and $r$ have the same word-label, $p$ has a large $\mathcal{P}$--component $w$. }
\label{fig:01}
\end{figure}
By the Bounded Coset Penetration property,
$r$ has a $P_i$--component $v$ of $S$--length at least
\[D-2|f|_S-c'-2\epsilon (1,0, |f|_S) = D-c.\]
Since $p=f^{-1}r$ and $r$ have the
same labels, it follows that $p$ has a $P_i$--component $w$ of $S$--length
at least $D-c$.

\end{proof}

\begin{corollary}\cite[Lemma 5.4]{MP07} \label{i:ip} 
Suppose that $G$, $Q$, $P$, and $R$ are as in the hypothesis of
Theorem~\ref{thm:2} and that $Q\cap P$ is a proper subgroup of $R$.

If  $\{K_1 , \dots , K_n \}$ is a collection of representatives of the maximal infinite parabolic subgroups of $Q$ up to conjugacy in $Q$ so that $K_1 =P\cap Q$, 
then $\{R, K_2, \dots , K_n \}$ is a collection of representatives of the maximal parabolic subgroups of $H$ up to conjugacy in $H$.
\end{corollary}
\begin{proof}
By Theorem~\ref{thm:2}~\eqref{i:pc}, a maximal parabolic subgroup of $H$  is conjugate to $R$ or $K_i$ for some $i=2,\dots n$. 
Hence $\{R, K_2, \dots , K_n \}$ is a collection of representatives of maximal parabolic subgroups.  That all these subgroups are different up to conjugacy follows
from the algebraic structure of $H$ as an amalgamated product.  In particular, since $K_i$ and $K_j$ are not conjugate in $Q$, they are not conjugate in $Q\ast_{Q\cap R} R$.  The subgroup $R < Q\ast_{Q\cap R} R$ is not conjugate to a subgroup of $Q$ since $Q\cap R$ is a proper subgroup of $R$.
\end{proof}

\begin{proof}[Proof of Theorem~\ref{thm:fully-quasiconvex}]

Suppose that every subgroup in $\mathcal{P}$ is LERF and slender.  Let $Q$ be a relatively quasiconvex subgroup of $G$ and let $g$ be an
element of $G$ not in $Q$.  Let $\{K_1, \dots , K_n\}$ be a collection of representatives of maximal infinite parabolic subgroups of $Q$ up
to conjugacy in $Q$; such a collection exists by Theorem~\ref{prop:ParabolicClasses}.  By Proposition~\ref{prop:parabolic-subgroups}, for each $K_i$ there
is a peripheral subgroup $P_i \in \mathcal{P}$ and $f_i \in G$ such that $K_i < P_i^{f_i}$.

We will construct an ascending sequence of relatively quasiconvex subgroups 
\[Q=Q_0 < Q_1 < \dots < Q_n=H\]
such that for each $k\in \{1, \dots , n\}$ the following properties hold.
\begin{enumerate}
\item\label{p-1} For $0 \leq j \leq k$,  the subgroup $Q_k \cap P_j^{f_j}$ is finite index in $P_j^{f_j}$.
\item\label{p-2} $\{ Q_k \cap P_1^{f_1}, \dots , Q_k \cap P_k^{f_k},  K_{k+1}, \dots , K_n \}$ is a collection of representatives of the maximal parabolic subgroups of $Q_k$ up to conjugation in $Q_k$.
\item\label{p-3} $g \not \in Q_k$.
\end{enumerate}
It will follow that $H=Q_n$ is a fully quasiconvex subgroup which
contains $Q$ and does not contain the element $g$.

Choose a geodesic $p$ from $1$ to $g$ in the relative Cayley graph
$\Gamma(G, \mathcal{P}, S)$.  Let $L$ be an upper bound for the
$S$--length of the $P_k$--components of the path $p$.

We now show how to construct $Q_k$, assuming that $Q_{k-1}$ has
already been constructed.  Let $C$ and $c$ be the constants provided by
Theorem~\ref{thm:2} for the subgroups $Q_{k-1}$ and $P_k^{f_k}$.
Since $P_k^{f_k}$ is slender, $K_{k} = Q_{k-1}\cap P_k^{f_k}$ is
finitely generated.  Define a finite set $F\subset
P_k^{f_k}\setminus K_k$ by 
\[ F  = \left\{ 
\begin{array}{ll}
  \{p\in P_k^{f_k}\setminus K_k\mid |p|_S\leq L+C+c\}\cup \{g\} &
  \mbox{if }g\in P_k^{f_k}\\
  \{p\in P_k^{f_k}\setminus K_k\mid |p|_S\leq L+C+c\} & \mbox{otherwise.}
\end{array}
\right.\]
Because $P_k^{f_k}\cong P_k$ is LERF, we may find a
finite index subgroup $R_k$ of $P_k^{f_k}$ satisfying 
\begin{itemize}
\item $K_k<R_k$, and
\item $f\notin R_k$ for all $f\in F$.
\end{itemize}
In particular, $|h|_S>L+C+c$ for any $h\in R_k\setminus Q_{k-1}$.
Let $Q_k=\langle Q_{k-1} \cup R_{k} \rangle$.  Note that the
hypotheses of the combination
Theorem \ref{thm:2} (and hence those of Corollary \ref{i:ip})
are satisfied for the relatively quasiconvex subgroup $Q_{k-1}$ and the parabolic subgroup
$R_k$.

We now verify properties \eqref{p-1}--\eqref{p-3} for the subgroup
$Q_k$ just constructed.  Property \eqref{p-1} follows from the fact
that $Q_k$ contains $Q_j \cap P_j^{f_j}$ for each $j$ between $1$ and $k-1$, and also contains $R_k$.
Corollary~\ref{i:ip} implies that
property~\eqref{p-2} is satisfied.  By Theorem~\ref{thm:2}\eqref{i:ge}
any geodesic in the Cayley graph $\Gamma(G, S\cup \bigcup
\mathcal{P})$ from $1$ to an element of element of $Q_k$ which is not
in $Q_{k-1}\cup R_k$ has a $P_k$--component of $S$--length greater
than $L+C$; it follows that the element $g$ does not belong to $Q_k$,
and so property \eqref{p-3} is also satisfied. 
This concludes the construction of the group $Q_k$, and the theorem
follows by taking $H=Q_n$.
\end{proof}

\section{Fillings of Relatively Hyperbolic Groups}\label{sec:fill}

Let $G$ be torsion free and relatively hyperbolic, relative to a
collection of subgroups $\mathcal{P} = \{P_1, \dots , P_m\}$, and let
$S\subset G$ be a finite generating set of $G$. 
Suppose that $S\cap P_i$ is a generating set of $P_i$ for each $i$.  


\begin{definition}
A \emph{filling}  of $G$ is determined by a collection of subgroups $\{N_i\}_{i=1}^m$ such that for each $i$,
$N_i$ is a normal subgroup of $P_i$; these subgroups are called \emph{filling kernels}. 
The quotient of $G$ by the normal subgroup generated by $\bigcup_{i=1}^m N_i$ is denoted by $G(N_1, \dots , N_m)$.
\end{definition}

The following result is due (in the present setting) independently to
Groves and Manning \cite[Theorem 7.2 and Corollary 9.7]{GM06} and to
Osin \cite[Theorem 1.1]{Os06-1}.  (Osin
actually proves a more general result, in which $G$ may have
torsion.)

\begin{theorem} \cite{GM06, Os06-1} \label{thm:4}
Let $F$ be a finite subset of $G$.  There exists a constant $B$
depending on $G$, $\mathcal{P}$, $S$, and $F$ with the following
property. 
If a collection of filling kernels $\{N_i\}_{i=1}^m$
satisfies $|f|_S>B$ for every nontrivial $f\in \cup_i N_i$, then
\begin{enumerate}
\item the natural map $\imath_i : P_i/N_i \longrightarrow G(N_1, \dots ,N_m)$ is injective,
\item\label{4hypq} $G(N_1, \dots ,N_m)$ is relatively hyperbolic relative to $\{
  \imath_i (P_i / N_i) \}_{i=1}^m$, and
\item\label{4inj} the projection $G \longrightarrow G(N_1, \dots ,N_m)$ is injective on $F$.
\end{enumerate}  
\end{theorem}

\subsection{Fillings and Quasiconvex Subgroups}
Let $H$ be a relatively quasiconvex subgroup of $G$.  A filling $G
\longrightarrow G(N_1, \dots , N_m)$ is an \emph{$H$--filling} if
whenever $H \cap P_i^g$ is non-trivial, $N_i^g \subset P_i^g \cap H$.

\begin{theorem} \cite[Propositions 4.3 and 4.5]{AGM08} \label{thm:3}
Let $H<G$ be a finitely generated relatively quasiconvex subgroup and $g \in G\setminus H$.
There is a finite subset $F\subset G$ depending on $H$ and $g$ with the following property.

If $\pi\co G \longrightarrow G(N_1, \dots , N_m)$ is an $H$--filling
which is injective on $F$, then $\pi (H)$ is a relatively quasiconvex subgroup of
$G(N_1, \dots , N_m)$, and $\pi (g) \not \in \pi (H)$.
\end{theorem}

\begin{remark}
  Propositions 4.3 and 4.5 in \cite{AGM08} use a different definition of relative quasiconvexity than the one we use here.
  In particular, their definition requires the subgroup to be finitely generated. We show in Appendix \ref{app:qcequiv} (specifically Corollary
  \ref{c:qcequiv}), that the definition used in \cite{AGM08} is equivalent to the one we use here under the assumption that the subgroup is finitely generated.
\end{remark}

\begin{proof}[{Proof of Theorem~\ref{thm:1}}]

Suppose that every subgroup in $\mathcal{P}$ is LERF and slender.  Let
$Q$ be a relatively quasiconvex subgroup of $G$ and let $g$ be an
element of $G$ not in $Q$.  By Theorem~\ref{thm:fully-quasiconvex},
there is a fully quasiconvex subgroup $H$ which contains $Q$ and does
not contain $g$.  We must choose filling kernels
$\{N_i\}_{i=1}^m$ such that $\pi\co G \longrightarrow G(N_1, \dots ,
N_m)$ is an $H$--filling, $G(N_1, \dots , N_m)$ is a word hyperbolic
group, $\pi (H)$ is a quasiconvex subgroup, and $\pi (g) \not \in \pi
(H)$.

By Theorem~\ref{prop:ParabolicClasses}, there is a collection
$\{K_1, \dots , K_n\}$ of representatives of the infinite maximal
parabolic subgroups of $H$ up to conjugacy in $H$.  
For each $r \in \{1, \dots , n\}$ there is an integer 
$i_r \in \{1, \dots , m\}$ 
and an element $f_r \in G$ such that $K_r^{f_r}$ is a finite index
subgroup of $P_{i_r}$.  The index $i_r$ is determined by $r$, but
there may be many distinct $f_r$ with this property.  On the other
hand, there will be only finitely many conjugates ${K_i}^g$ in
$P_{i_r}$.  Let $I_r$ be the intersection of these conjugates;  the
group $I_r$ is a finite index normal subgroup of $P_{i_r}$.
For $i \in \{1, \dots , m\}$ define the subgroup $M_i$ of $P_i$ as
\begin{displaymath}
M_i = \left \{ 
\begin{array}{ll}
P_i 												& 		\textrm{if $\{r \mid i_r = i \} = \emptyset $ } \\
													&																	\\
\bigcap \{ I_r \mid  i_r=i \} 	& 		\textrm{if $\{r \mid i_r = i \} \not = \emptyset $ }.
\end{array} 
\right .
\end{displaymath}
Put another way, if some conjugate of $H$ intersects $P_i$ nontrivially,
\[M_i = \bigcap\{{K_{r}}^g\mid {K_r}^g\cap {P_i}\neq\{1\}\mbox{, }r\in\{1,\ldots,n\}\mbox{, and }g\in G\}.\]
In other words, $M_i$ is the intersection of $P_i$ with all the
${K_r}^g$ which lie in $P_i$.
Because $M_i$ is a finite intersection of finite index normal subgroups 
of $P_i$, $M_i$ is a finite index normal subgroup of $P_i$.
From the definition, $M_i<{K_r}^g$ whenever $r\in \{r \mid i_r = i \}$
and ${K_r}^g<P_i$.

By Theorem~\ref{thm:3}, there is a finite subset $F \subset G$ such
that if $\pi\co  G\longrightarrow \bar{G}$ is any $H$--filling
which is injective on $F$, then $\pi (g) \not \in \pi (H)$ and $\pi
(H)$ is relatively quasiconvex.
Let $B>0$ be the constant from Theorem \ref{thm:4}, applied to this
finite subset $F$.

Since each $P_i$ is residually finite, there is some finite index
$\hat{N_i}\lhd P_i$ so that  $|p|_S>B$ for all $p\in
\hat{N}_i\setminus\{1\}$.  Let $N_i = \hat{N}_i\cap M_i$.  Since $N_i$
is an intersection of two finite index normal subgroups of $P_i$,
$N_i$ is a finite index normal subgroup of $P_i$.  Moreover, $|p|_S>B$
for every $p\in N_i\setminus\{1\}$.

\begin{lemma}\label{Hfill}
$G\to G(N_1,\ldots,N_m)$ is an $H$--filling.
\end{lemma}
\begin{proof}
Let $P_i\in \mathcal{P}$, $f\in G$, and suppose that $H\cap {P_i}^f\neq
\{1\}$.  
We must show that ${N_i}^f< H$.  Obviously it suffices to show that
${M_i}^f<H$. 
Since $G$ is torsion free, $H\cap {P_i}^f$ is infinite and
$H\cap{P_i}^f = {K_r}^{h}$ for some $h\in H$.  The group 
\[{M_i}^f = \bigcap\{{K_{r}}^g\mid {K_r}^g\cap {P_i}^f\neq\{1\}\mbox{,
}r\in\{1,\ldots,n\}\mbox{, and }g\in G\}\]
must therefore be contained in $H\cap {P_i}^f$, and thus in $H$.
\end{proof}

We now claim the subgroup $H<G$ and the
filling $G\to\bar{G} := G(N_1,\ldots,N_m)$ satisfy the
conclusions of Theorem \ref{thm:1}.  Indeed, conclusion \eqref{1qinh}
of Theorem \ref{thm:1}
is satisfied by construction.  By Theorem \ref{thm:4}.\eqref{4hypq},
the quotient $\bar{G}$ is  hyperbolic relative to a collection of
finite groups, hence $\bar{G}$ is hyperbolic.  Conclusion \eqref{1wh}
is thus established.
According to Lemma
\ref{Hfill}, $G\to\bar{G}$ is an $H$--filling. 
By Theorem~\ref{thm:4}.\eqref{4inj}, $G\to \bar{G}$ is injective on $F$.
Since the peripheral subgroups are slender, $H$ is finitely generated~\cite[Corollary 9.2]{Hr08}. 
We therefore may apply Theorem \ref{thm:3} to obtain conclusions \eqref{1qc} and
\eqref{1sep} of Theorem \ref{thm:1}.  Having established all the
conclusions, we have proved the theorem.
\end{proof}

\section{Applications to $3$--manifolds}\label{sec:3mfd}
\subsection{Kleinian groups}
Recall that a \emph{Kleinian group} is a discrete subgroup of
$\mathrm{Isom}(\mathbb{H}^3)$, the group of isometries of hyperbolic
$3$--space.  In this subsection we show that if all hyperbolic groups
are residually finite, then all finitely generated Kleinian groups are LERF\@.

\begin{proof}[Proof of Corollary \ref{cor:Kleinian}]
By Selberg's lemma \cite{Al}, every finitely generated Kleinian group
contains a torsion-free subgroup of finite-index.  Applying Corollary
\ref{l:finindex} it suffices to
consider torsion-free Kleinian groups.  Let $G$ be a torsion-free
Kleinian group and let $H<G$ be a finitely generated subgroup.
There are two cases:
\setcounter{case}{0}
\begin{case}
  $\HH^3/G$ has finite volume.
\end{case}
Either $H$ is geometrically finite or not.  Suppose first that $H$ is
not geometrically finite.  By a result of Canary \cite[Corollary
  8.3]{Ca96} together with the positive solution to the Tameness
Conjecture \cite{Ag,CG}, $H$ must be a \emph{virtual fiber subgroup}
of $G$.  This implies that there is a finite index subgroup of $G_0<G$
whose intersection $H_0$ with $H$ is normal in $G_0$, and so that
$G_0/H_0\cong \mathbb{Z}$.  The group $H_0$ is obviously separable in
$G_0$, so $H$ is separable in $G$ by Corollary \ref{l:finindex}.

If $H$ is geometrically finite, then it is relatively quasiconvex,
by Theorem \ref{thm:HK}, and we may apply Theorem \ref{t:c1}.
\begin{case}
  $\HH^3/G$ has infinite volume.
\end{case}
In this case it follows from the Scott core theorem and 
Thurston's geometrization theorem for Haken manifolds that  
$G$ is isomorphic to a
geometrically finite Kleinian group $G'$ (see, for example,
\cite[Theorem 4.10]{mt:hmkg}). 
If $H'$ is the image of $H$ in
$G'$ then obviously $H$ is separable in $G$ if and only if $H'$
is separable in $G'$.  We therefore may as well assume that $G$
is geometrically finite to begin with.

Since $G$ is geometrically finite and infinite covolume, every
finitely generated subgroup of $G$ is geometrically finite by an
argument of Thurston (see \cite[Proposition 7.1]{morgan:uniform} or
\cite[Theorem 3.11]{mt:hmkg} for
a proof).
In particular, $H$ is geometrically finite, so Theorem \ref{thm:HK}
implies that $H$ is a relatively quasiconvex subgroup of $G$.  
To finish, we apply Theorem \ref{t:c1} again.
\setcounter{case}{0}
\end{proof}

\subsection{Subgroup separability in finite volume $3$--manifolds}
Here we prove that if compact hyperbolic $3$--manifold groups are QCERF
then finite volume hyperbolic $3$--manifold groups are LERF\@.  
\begin{proposition}\label{c:fv}
  If all fundamental groups of compact hyperbolic $3$--manifolds are
  QCERF, then all fundamental groups of finite volume hyperbolic $3$--orbifolds
  are LERF\@.
\end{proposition}
\begin{proof}
Let $G = \pi_1(M)$, where $M$ is some finite volume hyperbolic $3$--orbifold.
Applying Corollary \ref{l:finindex},
we may pass to a finite cover and
assume that $M$ is an orientable manifold.
It follows that $G$ is relatively hyperbolic, relative to some
finite collection $\mathcal{P}=\{P_i,\ldots,P_m\}$ of subgroups, each of which is isomorphic
to $\ZZ\oplus\ZZ$.

Now suppose $Q<G$ is some finitely generated subgroup.
We must show that $Q$ is separable.
If $Q$ is geometrically infinite, then we argue as we did in the proof
of \ref{cor:Kleinian} that $Q$ is the
fundamental group of a virtual fiber, and thus separable.  We may thus
suppose that $Q$ is geometrically finite, and therefore a relatively quasiconvex
subgroup of $G$.  

Let $g\in G\setminus Q$.
We then apply Theorem \ref{thm:fully-quasiconvex} to enlarge $Q$ to a
fully quasiconvex subgroup $H$ not containing $g$.  Let $F\subset G$ be the
finite set obtained by applying Theorem \ref{thm:3} to $G$, $H$, and
$g$.  Let $B_1$ be the constant from Theorem \ref{thm:4} applied to $G$
and $F$, with respect to some generating set $S$ for $G$.  We will
choose a cyclic filling kernel $N_i<P_i$ for each $P_i\in \mathcal{P}$.
The Hyperbolic Dehn Surgery Theorem of Thurston \cite{Th,HK} implies
there is some constant $B_2$ so that if the generators of the $N_i$
are chosen to have length greater that $B_2$, then the filling
$G(N_1,\ldots,N_m)$ will be the (orbifold) fundamental group of a
hyperbolic orbifold obtained by attaching orbifold solid tori to the boundary
components of a compact core of $M$.
Let $B = \max\{B_1,B_2\}$.

For each
$P_i\in \mathcal{P}$ we choose 
some cyclic $N_i<P_i$.
For each $i$ let $n_i\in P_i$ satisfy $|n_i|_S>B$.
There are at most finitely many conjugates
${P_i}^{t_1},\ldots {P_i}^{t_k}$ so that 
${P_i}^{t_j}\cap H$ is
nonempty; for each such $j$, the group ${P_i}^{t_j}\cap H$ is finite index
in ${P_i}^{t_j}$.   If there are no such conjugates, we choose $N_i = \langle
n_i \rangle$.
Otherwise, we let $N_i = \langle n_i^\alpha \rangle$, where the power
$\alpha\in \mathbb{N}$
is chosen so that $n_i^\alpha\in H^{{t_j}^{-1}}$ for each $j$, and so $|n_i^\alpha|_S>B$.

With the $N_i<P_i$ chosen as above, the $H$--filling
$G(N_1,\ldots,N_m)$ is a compact hyperbolic $3$--orbifold group.  Let
$\pi\co G\to G(N_1,\ldots,N_m)$ be the quotient map.  By Theorem
\ref{thm:3}, $\pi(H)$ is a quasiconvex subgroup of
$G(N_1,\ldots,N_m)$, not containing $\pi(g)$.  By assumption, compact
hyperbolic $3$--manifold groups (and thus compact hyperbolic
$3$--orbifold groups) are QCERF\@.  There is therefore some finite group
$F$, and some $\phi\co G(N_1,\ldots,N_m)\to F$ with
$\phi(\pi(g))\notin\phi(\pi(H))$.  Since $\phi(\pi(H))$ contains
$\phi(\pi(Q))$, we have separated $g$ from $Q$ in a finite quotient.
\end{proof}

\begin{remark}
  No part of the proof of Proposition \ref{c:fv} rests in an essential way on the
  results in this paper or in \cite{AGM08}, but only on facts about
  hyperbolic $3$--manifolds and $3$--orbifolds which could be deduced
  by geometric arguments, based on the Gromov-Thurston $2\pi$ Theorem.
  On the other hand, Proposition \ref{c:fv} is a nice illustration of
  the general principle that a relatively hyperbolic group which can
  be approximated by QCERF Dehn fillings must itself be QCERF\@.  
\end{remark}

\section{Acknowledgments}
Thanks very much to Ian Agol and Daniel Groves, both for useful
conversations, and for helpful comments on an earlier draft of this
note.  We also thank the referee for insightful comments and corrections.

Manning was partially supported by NSF grant DMS-0804369. 
Mart\'inez-Pedroza was funded as a Britton Postdoctoral Fellow at
McMaster University, grant A4000.

\appendix
\section{On the equivalence of definitions of quasiconvexity} \label{app:qcequiv}

The current paper relies heavily on results about relatively
quasiconvex subgroups of relatively hyperbolic groups proved in the
papers \cite{MP07} and \cite{AGM08}.  These papers unfortunately use
different definitions of relative quasiconvexity, but we show in this
appendix that the two definitions agree, at least for finitely
generated subgroups.

First, a few words about the literature.
Dahmani \cite{Da03} and Osin \cite{Os06} studied classes of
subgroups of relatively hyperbolic groups which they called
\emph{relatively quasiconvex}, intending to generalize the notion of
quasiconvexity in hyperbolic groups.  Dahmani's definition was a
dynamical one, whereas Osin's was Definition \ref{def:O};  Osin's
definition was used in \cite{MP07}.
Hruska in ~\cite{Hr08} gave several definitions of relative
quasiconvexity in the setting of countable (not necessarily finitely
generated) relatively hyperbolic groups, including
definitions based on Osin's and Dahmani's, and showed that they are
equivalent.
The authors of \cite{AGM08} were mainly interested in relatively
hyperbolic structures on groups which were already hyperbolic, and
used a definition of relative quasiconvexity (based more closely on
the usual metric notion of quasiconvexity) different from any of those
in \cite{Hr08}.  The definition in \cite{AGM08} applies only to a
finitely generated subgroup of a finitely generated relatively
hyperbolic group.

Throughout this section $G$ will be relatively hyperbolic, relative to a finite collection of subgroups 
$\mc{P} = \{ P_1,\dots ,P_n\}$,  and $S$ will be a finite generating
set for $G$.  
The \emph{cusped space} (recalled below) for $G$ with respect to
$\mc{P}$ and $S$ will be denoted by $X(G,\mc{P},S)$, 
and $d( \cdot , \cdot )$ will denote the path metric on the cusped
space.
In particular, we will not 
need the word metric on $G$, but only the metric induced by this path metric.
For more detailed definitions and background on
cusped spaces for relatively hyperbolic 
groups we refer the reader to~\cite{GM06} and~\cite{Hr08};  we sketch
the construction and recall some terminology here for the reader's
convenience. 

Let $A$ be a discrete metric space with metric $\rho$.  The \emph{combinatorial horoball based
  on $A$} is a graph $\mc{H}(A)$ with vertex set $A\times\ZZ_{\geq 0}$, so that
\begin{itemize}
\item $(a,n)$ is connected by an edge to $(a,n+1)$ for any $a\in A$,
  and
\item for $n\geq 1$, $(a,n)$ is connected to $(a',n)$ whenever
  $\rho(a,a')\leq 2^n$.
\end{itemize}
Edges of the first type are called \emph{vertical}; edges of the
second type are called \emph{horizontal}.
We say that a vertex $(a,n)$ of $\mc{H}(A)$ has \emph{depth $n$}.
If $A$ is a subset of a path metric space $Y$, we may \emph{attach a
  horoball to $Y$ along $A$} by gluing $A\subseteq Y$ to $A\times
\{0\}\subset \mc{H}(A)$, and taking the obvious path metric on the
union.  If $G$ is finitely generated by $S$, we take $Y$ to be the
Cayley graph of $G$ with respect to $S$; any subset of $G$ inherits a
discrete metric from the path metric on the Cayley graph.
The \emph{cusped space}
$X(G,\mc{P},S)$ is the space obtained by simultaneously attaching
horoballs to $Y$ along all left cosets $tP$ for $P\in \mc{P}$.  A 
vertex $v$ of $X(G,\mc{P},S)$ corresponds either to a group element, if
it lies in the Cayley graph of $G$, or otherwise to a triple
$(tP,g,n)$ where $tP$ is a left coset of an element of $\mc{P}$, the
element $g$ lies in $tP$, and $n>0$ is the depth of $v$ in the
attached horoball $\mathcal{H}(tP)$.  In what follows we
do not distinguish between $v$ and the corresponding group element or
triple.

The group $G$ is relatively hyperbolic, relative to $\mc{P}$, if and
only if the space $X(G,\mc{P},S)$ is Gromov hyperbolic (see for example
\cite[Theorem 3.25]{GM06}).  If so, then
$G$ acts on $X(G,\mc{P},S)$ \emph{geometrically
  finitely},\footnote{Hruska uses the term ``cofinitely.''}
  meaning in particular:
\begin{enumerate}
\item  Given any $n\geq 0$, let $X_n$ be the subset of $X(G,\mc{P},S)$
  obtained by deleting all vertices of height greater than $n$.  Then
  $G$ acts cocompactly on $X_n$ (which is an example of what Hruska calls a \emph{truncated
    space} for the action of $G$ on $X(G,\mc{P},S)$).
\item  For fixed $n$, there are only finitely many components of  of
  $X(G,\mc{P},S)\setminus X_n$, up to the action of $G$.  (The
  components of $X(G,\mc{P},S)\setminus X_{n-1}$ are called
  \emph{$n$--horoballs}, for $n\geq 1$.  
  If $n$ is understood, we call them \emph{horoballs}.  A
  $0$--horoball is a $1$--neighborhood of a $1$--horoball, and is
  equal to $\mc{H}(tP)$ for some coset $tP$ of some $P\in
  \mc{P}$.)
\end{enumerate}

\subsection*{Relatively Quasiconvex Subgroups according to
  Agol--Groves--Manning.}
Suppose that $H$ is a relatively hyperbolic group, and let $\mc{D} = \{ D_1, \dots , D_m \}$ be the
peripheral subgroups of $H$ and $T$ a finite generating set for $H$.  Let $\phi\co H\to G$ be a homomorphism.  If every
$\phi(D_i)\in \mc{D}$ is conjugate in $G$ into some $P_j\in \mc{P}$, we say that the map $\phi$ \emph{respects the peripheral structure on} $H$.

Given such a map $\phi$, one can extend it to a map $\check{\phi}$
between zero-skeletons of cusped spaces in the following way:  
For each $D_i \in \mc{D}$, choose an element $c_i \in G$ (of minimal length) and some $P_{j_i}$ such that $\phi(D_i) \subseteq c P_{j_i} c^{-1}$. 
For $h\in H$, $\check\phi(h) = \phi(h)$. 
For a vertex $(sD_i,h,n)$ in a horoball of $X(H,\mc{D},T)$, define
\[ \check\phi(sD_i,h,n) = (\phi(s)c_iP_{j_i}, \phi(h)c_i, n) .\]

\begin{lemma}\label{lem:AGM}\cite[Lemma 3.1]{AGM08}
Let $\phi\co H\to G$ be a homomorphism which respects the peripheral structure on $H$.
The extension $\check{\phi}$ defined above
is $H$--equivariant and lipschitz.
If $\phi$ is injective, then $\check{\phi}$ is proper.
\end{lemma}

Recall that, for $C\geq 0$,
a subset $A$ of a geodesic metric space $X$ is \emph{$C$--quasiconvex} if
every geodesic with endpoints in $A$ lies in a $C$--neighborhood of
$A$.  The subset is \emph{quasiconvex} if it is $C$--quasiconvex for
some $C$.
The following is a slight paraphrase of the definition from
\cite{AGM08}:
\begin{definition} (QC-AGM) \cite[Definition 3.11]{AGM08} \label{def:AGM}
Let $G$ be as above, and let $H<G$ be finitely generated by a set $T$.
We say that $H$ is 
\emph{(QC-AGM) relatively quasiconvex} in
$(G,\mathcal{P})$ if, for some finite collection of subgroups $\mc{D}$
of $H$,
\begin{enumerate}
\item $H$ is relatively hyperbolic, relative to $\mc{D}$, and
\item if $\iota\co H\to G$ is the inclusion, then the map $\check\iota\co
  X(H,\mc{D},T)^0\to X(G,\mathcal{P},S)^0$ described in Lemma
  \ref{lem:AGM} has quasiconvex image.
\end{enumerate}
\end{definition}

\subsection*{Relatively Quasiconvex Subgroups according to Hruska.}
The following definition is direct from Hruska \cite{Hr08}, where it
is called \emph{QC-3}.  Hruska shows in \cite{Hr08} that this
definition is equivalent to several others, including our
Definition \ref{def:O}. 
\begin{definition}(QC-H)\label{def:H}\cite[Definition 6.6]{Hr08}
A subgroup $H \le G$ is \emph{(QC-H) relatively quasiconvex} if the following holds.
Let $(X,\rho)$ be some (any) proper Gromov hyperbolic space on
which $(G,\mc{P})$ acts geometrically finitely.
Let $X-U$ be some (any) truncated space for $G$ acting on $X$.
For some (any) basepoint $x \in X - U$
there is a constant $\mu\ge 0$ such that whenever $c$ is a geodesic in $X$
with endpoints in the orbit $Hx$, we have
\[
   c \cap (X-U) \subseteq \nbd{Hx}{\mu},
\]
where the neighborhood is taken with respect to the metric $\rho$ on $X$.
\end{definition}
\begin{remark}\label{r:nondef}
  The meaning of ``some (any)'' in Definition \ref{def:H} just means
  that the word ``some'' can be replaced by ``any'' 
  without affecting which subgroups of $G$ are (QC-H) relatively
  quasiconvex.  Thus ``Definition'' \ref{def:H} has some
  non-definitional content, established in \cite[Proposition 7.5 and
  7.6]{Hr08}.)
\end{remark}
\begin{definition}
  Let $A\subset X=X(G,\mc{P},S)$ be a horoball, and let $R>0$.
 We say that a geodesic \emph{$\gamma$ penetrates the horoball $A$ to depth $R$} 
 if there is a point 
$p \in \gamma\cap A$ at distance at least $R$ from $X\setminus A$. 
We say
  that \emph{$A$ is $R$--penetrated by the subgroup $H$}
  if there is a geodesic $\gamma$
  with endpoints in $H$ penetrates the horoball $A$ to depth $R$.
\end{definition}
The goal of this subsection is to prove the following proposition.

\begin{proposition}\label{l:shallowtorsion}
Let $H<G$ be (QC-H) relatively quasiconvex.
Then there is a constant $R=R(G, \mc{P}, S, H)$, so that whenever 
a $0$--horoball is $R$--penetrated by $H$, the intersection of $H$ with
the stabilizer of that horoball is infinite.
\end{proposition}

Before  the proof, we quote a proposition from \cite{Hr08} and prove two
lemmas.

\begin{proposition}\cite[Proposition 9.4]{Hr08}\label{prop:cosets}
\label{prop:CloseCosets}
Let $G$ have a proper, left invariant metric $d$, and suppose $xH$ and $yK$ are arbitrary left cosets of subgroups of $G$.
For each constant $L$ there is a constant $L'=L'(G,d,xH,yK)$ so that in the metric space $(G,d)$ we have
\[
   \nbd{x H}{L} \cap \nbd{y K}{L} \subseteq   \nbd{x H x^{-1} \cap y K y^{-1}}{L'}.
\]
\end{proposition}

\begin{lemma}\label{lem:H-1}
Let $H$ be a (QC-H) relatively quasiconvex subgroup of $G$.  Let $A$
be a $0$--horoball of $X(G,\mc{P},S)$, whose stabilizer is $P^t$ for $P\in
\mc{P}$.  If $A$ is $R$--penetrated by $H$ for all $R>0$, then $H\cap
P^t$ is infinite.
\end{lemma}

\begin{proof}
It suffices to show that, for every $M >0$, there is some $h$
in $H\cap P^t$ with $d(1,h)>M $.

Let $\mu$ be the
quasiconvexity constant of Definition~\ref{def:H}
for $H$ and the space $X'$ which consists of all vertices in $X(G, \mc{P}, S)$ at depth $0$.  Let $C$ be the constant given by Proposition~\ref{prop:cosets}  such that 
\[ \nbd{H}{\mu} \cap  tP  \subseteq  \nbd{H\cap tPt^{-1}}{C} ,\]
where the neighborhoods are taken in the cusped space.

Suppose that $\gamma$ is a geodesic with endpoints in $H$ which
penetrates the horoball $A$ to depth $M +C$.  The first and last points
of $\gamma\cap A$ are group elements, $a$ and $b$, both in the coset $tP$.
Since $H$ is (QC-H) relatively quasiconvex, $a$ and $b$ are elements
of $\nbd{H}{\mu} \cap  tP$ and therefore (using
Proposition~\ref{prop:cosets}) 
there are elements $h_1$ and $h_2$ 
in $H\cap P^t$  such that 
$d(h_1,a) \leq C$ and $d(h_2, b) \leq C$.
Since $d(a, b)\geq 2(M +C)$, 
\[d(1,
h_1^{-1}h_2) = d(h_1, h_2) \geq  2(M +C) - 2C \geq 2M >M .\qedhere\] 
\end{proof}

\begin{lemma}\label{lem:H-2}
Let $H$ be a (QC-H) relatively quasiconvex subgroup of $G$. 
Let $\mu$ be the quasiconvexity constant of Definition~\ref{def:H} 
for $H$ and the space $X'=X_0$ which is obtained from $X(G,\mc{P},S)$
by deleting all vertices of positive depth.

Let $R>0$, and
let $A$ be a $0$--horoball, stabilized by $P^t$ for $P\in \mc{P}$.
If  $A$ is $R$--penetrated by $H$, then 
there is a horoball $A'$ so that
\begin{enumerate}
\item $A' = hA$ for some $h\in H$,
\item $d(A',1)\leq \mu$, and
\item $A'$ is $R$--penetrated by $H$.
\end{enumerate}
\end{lemma}
\begin{proof}
Suppose that $\gamma$ is a geodesic with endpoints $h_1$ and $h_2$ in $H$ 
which penetrates the horoball $A$ to depth $R$.  Let $a$ and $b$ be
the first and last vertices of $\gamma\cap A$.
By (QC-H) relative quasiconvexity,
there is some $h\in H$ so that $d(a, h)\leq \mu$. 
The geodesic $h^{-1}\gamma$ goes between $h^{-1}h_2$ and $h^{-1}h_2$,
and penetrates the horoball $A'=h^{-1}A$ to depth $R$.  Moreover,
\[d(1,A')\leq d(1,h^{-1}a) = d(a,h) \leq \mu.\qedhere\]
\end{proof}

\begin{proof}[Proof of Proposition~\ref{l:shallowtorsion}]
Suppose there is no such number $R$. 
There must be a sequence of integers $R_i\to\infty$ and a
sequence of $0$--horoballs $\{A_i\}$, so that, for each $i$, the
horoball $A_i$ is $R_i$--penetrated by $H$, but the intersection of
the stabilizer of $A_i$ with $H$ is finite.

For $h\in H$, the stabilizer of $hA_i$ is conjugate (by $h$) to the
stabilizer of $A_i$.  Using Lemma~\ref{lem:H-2}, we can therefore
assume that $d(1,A_i)\leq \mu$ for each $i$.  By passing to a
subsequence, we can therefore assume that the sequence $\{A_i\}$ is
constant.  It follows that $A_0$ is $R_i$--penetrated by $H$ for all
$i$.  Lemma~\ref{lem:H-1} then implies that the intersection of $H$
with the stabilizer of $A_0$ is infinite, which is a contradiction.
\end{proof}

\subsection*{Equivalence of the two definitions.}
In this section, $G$ will be a relatively hyperbolic group, relative to a finite collection of subgroups $\mc{P}$, and $S$ will be a finite generating set for G.  Let   $X(G,\mc{P},S)$ be the cusped space for $G$ with respect to $\mc{P}$ and $S$, and let $\delta$ be its hyperbolicity constant.

\begin{theorem}
Let $H$ be a finitely generated subgroup of $G$.  Then $H$ is (QC-H)
relatively quasiconvex if and only if $H$ is (QC-AGM) relatively
quasiconvex.
\end{theorem}
\begin{proof}
One direction is easy.  Suppose that $H<G$ is (QC-AGM) relatively
quasiconvex, generated by the finite set $T$, and with peripheral
subgroups $\mc{D}$.
Recall that to define $ \check{\iota} \co X(H, \mc{D}, T)^0
\longrightarrow X(G, \mc{P}, S)^0$,  
an element $c_i\in G$ was chosen for each $D_i\in \mc{D}$ so that
$D_i \subset c_i P_{j_i} c_i^{-1}$.  Let 
$ C = \max\{ d(1,c_i) \mid D_i\in \mc{D}\}$, and
let $C_q$ be the constant of quasiconvexity in the
definition of (QC-AGM) quasiconvexity.  
As remarked at the beginning of the Appendix, the cusped
space $X = X(G,\mc{P},S)$ is acted on geometrically finitely by $G$, and
the subspace $X-U = X_0 \subset X(G,\mc{P},S)$ obtained by deleting
$1$--horoballs is a truncated space for the action.
Moreover, as explained in Remark \ref{r:nondef}, it suffices to find a
$\mu$ which works for this choice of $X$ and $X-U$, and for the $H$--orbit
of $1$ in $X$.  Let $x$, $y\in H$, and let $\gamma$ be any geodesic
joining them in $X$. 
Let $z$ be a vertex of $\gamma$ contained in $X_0$.
By (QC-AGM), there is some point $w\in \check\iota(X(H,\mc{D},T)^0)$ so
that
$ d(z,w)\leq C_q$.
It follows that $w\in X_{C_q}$, but any point in
$\check\iota(H,\mc{D},T)\cap X_{C_q} $
is at most $C + C_q$ away from some point in $H$.  It follows that $z$
is no further than $\mu:= C+2C_q$ from $H$, and so $H$ is (QC-H)
relatively quasiconvex.

We now establish the other direction.
Let $H$ be a subgroup of $G$, and suppose that $H$ is (QC-H)
relatively quasiconvex.  Let $\mc{D}$ be a collection of representatives of
the $H$--conjugacy classes of 
infinite maximal parabolic subgroups of $H$.  By~\cite[Theorem 9.1]{Hr08}, 
$H$ is relatively hyperbolic, 
relative to $\mc{D}$.
By Lemma~\ref{lem:AGM}, the inclusion $\iota\co  H\longrightarrow G$
extends to a lipschitz map of ($0$--skeletons of) cusped spaces
\[  \check{\iota} \co X(H, \mc{D}, T)^0 \longrightarrow X(G, \mc{P}, S)^0. \]
We need to prove that the image $Y = \check{\iota}( X(H,\mc{D}, T)^0)$
of $\check{\iota}$ is quasiconvex.  

Let $R$ be the constant provided by
Proposition~\ref{l:shallowtorsion} for the subgroup $H$.   
Let $X' = X_{100\delta +R}$ be the subspace of $X(G, \mc{P}, S)$ consisting of all vertices at depth at most $100\delta +R$. 
Since $H$ is (QC-H) relatively quasiconvex, there is a constant $\mu$ such that for any geodesic $\zeta$ in $X(G, \mc{P}, S)$ with endpoints in $H$,
\[ \zeta \cap X' \subset \nbd{H}{\mu} \subset \nbd{Y}{\mu},\]
where the neighborhoods are taken with respect to the metric on $X(G, \mc{P} ,S)$.

Let $x$ and $y$ be vertices of $Y$ and let $\gamma$ be a geodesic
between them.  We will show that the vertices of $\gamma$ are contained in the
$M$-neighborhood of $Y$, where $M$ is a constant independent of $x$,
$y$, and $\gamma$.  We divide the proof into five 
(not necessarily
disjoint) cases. 

\begin{case} \label{in.same.horoball}
The points $x$ and $y$ lie deeper than $10\delta$ in the same
horoball.
\end{case}
By recalling some easily verified properties of the geometry of horoballs,
we will show that $\gamma$ is contained in the $M_1$-neighborhood of $Y$, where
 \[ M_1= 6.\] 
To begin with, the
$10\delta$--horoball containing $x$ and $y$ is convex
(see~\cite[Lemma 3.26]{GM06}). 
Second, 
any geodesic with the same endpoints as $\gamma$ is Hausdorff distance
at most $4$ from $\gamma$.
Finally, there is a geodesic $\gamma'$ of a particularly nice form
with the same endpoints as $\gamma$.  The geodesic $\gamma'$ is a
\emph{regular geodesic}, which means that all its edges are vertical,
except for at most three consecutive horizontal edges at maximum depth
(see \cite[Lemma
3.10]{GM06}). 
Since the vertical subsegments of
$\gamma'$ start at points in $Y$ and are vertical, they stay in $Y$,
and so $\gamma'$ stays in a $2$--neighborhood of $Y$.  
As $\gamma$ is contained in a
$4$--neighborhood of $\gamma'$, we have $\gamma$ contained in a
$6$--neighborhood of $Y$.

\begin{case}\label{h.near.horoball}
The points $x$ and $y$ are elements of $H$, they are in the
neighborhood of radius $\mu$ of a horoball $\mc{H}(tP )$, and the 
geodesic $\gamma$ penetrates the horoball $\mc{H}(tP )$ to depth
larger than $100\delta + R$. 
\end{case}
In this case, we will approximate $\gamma$ by a regular geodesic
inside $\mc{H}(tP)$ with (possibly different) endpoints in $Y$.
Without loss of generality, assume that $x$ is the identity, and so
$d(1,t) \leq \mu$.  

By Proposition~\ref{l:shallowtorsion}, the intersection $H \cap P^t$
is infinite.  It follows that $H \cap P^t = D^s$ for some $D \in
\mathcal{D}$ and $s \in H$.  We claim $s$ can be chosen so that
$d(1,s)<K$ for a constant $K$ independent of $x$, $y$, and $\gamma$.
Indeed, we observe that
the set 
\[ W = \{ (r, P) \in G \times \mathcal{P}  \mid  d(1, r)\leq \mu,  \# ( H\cap P^r ) = \infty  \}\]
is finite.  For each  $w=(r, P)\in W$  choose $u_w
\in H$ so that $H \cap P^r = D^u$ for some $D\in \mc{D}$;
we let 
$K$ be the maximum of $d(1, u_w)$ over all $w\in W$. 

We further claim that there is an element $y' \in H \cap P^t$ such
that $d(y, y') \leq L$, for a constant $L$ independent of $x$, $y$,
and $\gamma$.  Indeed, for each $w = (r,P)$ in the set $W$ defined
above, Proposition~\ref{prop:cosets} implies we can find an $L_w>0$ so
that 
\[ H \cap \nbd{rP}{\mu} \subseteq \nbd{H \cap P^r}{L_w}; \]
we let $L$ be the maximum $L_w$ over all $w\in W$.

Recall that to define $ \check{\iota} \co X(H, \mc{D}, T)^0 \longrightarrow X(G, \mc{P}, S)^0$, 
an element $c_i\in G$ was chosen for each $D_i\in \mc{D}$ so that
$D_i \subset c_i P_{j_i} c_i^{-1}$.  Let 
\begin{equation}\label{eq:C} C = \max\{ d(1,c_i) \mid D_i\in \mc{D}\}.
\end{equation}
The subgroup $D$ is equal to $D_i$ for some $i$, and we set $c = c_i$
for the same $i$.

Consider the elements $(sD, s, 10\delta)$ and $(sD, y's, 10\delta)$
of $X(H, \mc{D}, T)$ and their corresponding images in $Y$ given by
$(scP, sc, 10\delta)$ and $(scP, y'sc, 10\delta)$.
The points $(scP, sc, 10\delta)$ and $(scP, y'sc, 10\delta)$ belong to
the same $10\delta$--horoball, which is convex in $X(G,\mc{P},S)$, as
we noted in Case \ref{in.same.horoball}.  Also as noted in Case
\ref{in.same.horoball}, there is a regular geodesic $\gamma'$ joining
the points  $(scP, sc, 10\delta)$ and $(scP, y'sc, 10\delta)$;  since
the endpoints lie in $Y$, the geodesic $\gamma'$ is contained in the
$2$--neighborhood of $Y$.

On the other hand, the endpoints of the geodesics $\gamma$ and
$\gamma'$ are close, namely, 
\[  d(1, (scP, sc, 10\delta)) \leq  d(1,s) + d(1, c) +10\delta \leq  K+ C +10\delta , \]
and
\[  d(y, (scP, y'sc, 10\delta)) \leq  d(y,y') + d(y',(scP, y'sc, 10\delta))  \leq   L + K + C +10\delta. \]
Since $X(G, \mc{P}, S)$ is $\delta$--hyperbolic,
the Hausdorff distance between $\gamma'$ and $\gamma$ is at most the
distance between endpoints plus $2\delta$.  Thus if 
\[ M_2 = 2 \delta +  K + L + C +10\delta, \]
then $\gamma$ is contained in the $M_2$-neighborhood of $Y$.
This completes this case.

\begin{case}\label{pf_case_3}
Suppose $x$ and $y$ are elements of $H$.
\end{case} 

We split $\gamma$ into subsegments $\gamma_1, \gamma_2, \dots
,\gamma_k$ such that no $\gamma_i$ contains any group element (depth
$0$ vertex) in its
interior, but the endpoints of each $\gamma_i$ are group elements.
Observe that each $\gamma_i$ is either a single edge or a
geodesic segment contained in a $0$--horoball.  Furthermore, since $H$ is (QC-H) relatively quasiconvex, the endpoints of each $\gamma_i$ are contained in the $\mu$-neighborhood of $H$. 
We claim that each $\gamma_i$ is contained in the $M_3$-neighborhood of $Y$, where
\[ M_3 = 110\delta + R + 3\mu + 2 + M_2 .\]
If $\gamma_i$ is an edge, then the claim is immediate, so we suppose
$\gamma_i$ is contained in a $0$--horoball $\mc{A}$.
First, suppose $\gamma_i$  does not penetrate $\mc{A}$
to
depth $110\delta + R + 2 \mu$. 
An easy argument shows that the length
of a geodesic in a combinatorial horoball is at most twice its maximum depth plus $4$, so we have
$|\gamma_i|< 220\delta+2R + 2\mu+4$, and $\gamma_i$ is therefore
contained in the 
$(110\delta + R + 3\mu+2)$--neighborhood of $H$. 
In particular, $\gamma_i$ is contained in $M_3$--neighborhood of $Y$.

Suppose on the other hand that $\gamma_i$ penetrates the horoball $\mc A$ to depth
$110\delta + R + 2\mu$.  Let $h_1$ and $h_2$ be elements of $H$ which
are at distance at most $\mu$ from the endpoints of $\gamma_i$, and
let $\alpha$ be a geodesic between them.  Since $X(G, \mc{P}, S)$ is
$\delta$--hyperbolic, the Hausdorff distance between $\gamma_i$ and
$\alpha$ is at most $2\delta + \mu$.  It follows that $\alpha$
penetrates the horoball $\mc A$ to depth $100\delta + R + \mu$, and
hence it satisfies the condition of Case~\ref{h.near.horoball}.
Therefore, $\gamma_i$ is in the $(2\delta + \mu + M_2)$--neighborhood
(and hence in the $M_3$--neighborhood) of $Y$. 

\begin{case}
Suppose $x$ and $y$ lie at depth no more than $50\delta$ in $X(G,\mc{P},S)$.
\end{case}

If $x\in Y$ lies in a $1$--horoball, then $x = (tP,hc_i,n)$ for some
$P\in \mc{P}$, some $h\in H$, some $i\in \{1,\ldots,m\}$, and some
$n\leq 50\delta$;
otherwise, $x\in H$.  In any case, there
is an element $h_1 \in H$ such that 
$d(x, h_1) \leq 50\delta + C$, where $C$ is the constant defined in \eqref{eq:C}.
By the same argument,
there is an element $h_2 \in H$ such that
$d(x, h_2) \leq 50\delta + C$.
Since $X(G, \mc{P}, S)$ is $\delta$--hyperbolic, the Hausdorff
distance between $\gamma$  and any geodesic $\gamma'$ between $h_1$ and $h_2$ is at most $52\delta +  C$.  
We may apply Case \ref{pf_case_3} to $\gamma'$, and deduce that
$\gamma$ is contained in the $M_4$--neighborhood of $Y$, where
\[ M_4 = 52\delta +  C + M_3. \] 

\begin{case}
Suppose either $x$ or $y$ lies inside a $50\delta$--horoball, but we are not in Case 1.
\end{case}

Here we follow the proof of the last case of ~\cite[Proposition 3.12]{AGM08}.
If $x$ or $y$ lies in a horoball, it is connected by a vertical path
to a point in the right coset $Hc_i$ at depth $0$ in $X(G, \mc{P},
S)$.   It is therefore possible to modify $\gamma$ (by appending and deleting
(mostly) vertical paths lying in a $3$--neighborhood of $Y$) to a
$10\delta$--local geodesic $\gamma'$ with endpoints within $C$ of
$H$; the geodesic $\gamma$ is contained in a $3$--neighborhood of
$\gamma'\cup Y$.
By~\cite[III.H.1.13(3)]{bridhaef}, $\gamma'$  is a 
$(\frac{7}{3}, 2 )$--quasi-geodesic.
Since quasi-geodesics track geodesics, there is a constant $L_Q$ depending only on $\delta$ and $C$, 
and a geodesic $\gamma''$ with endpoints in $H$ such that the Hausdorff distance between $\gamma'$ and $\gamma''$ is at most $L_Q$. 
By Case~\ref{pf_case_3}, $\gamma''$ is contained in the $M_2$-neighborhood of $Y$. Let 
\[ M_5 = 3 + L_Q +M_2, \]
and observe that $\gamma$ is contained in the $M_5$-neighborhood of $Y$.

Finally, we set $M = \max\{M_1,\ldots,M_5\}$, and note that $M$ does
not depend on the vertices $x$ and $y$ of $Y$, or on the geodesic
$\gamma$ joining them.  It follows that $Y=\check{\iota}( X(H,\mc{D},
T)^0)$ is $M$--quasiconvex in $X(H,\mc{D},T)$, and so $H$ is
(QC-AGM)
relatively
quasiconvex in $(G,\mc{P})$.
\setcounter{case}{0}
\end{proof}

Applying the main result of Hruska \cite{Hr08} on the equivalence of
various definitions of relative quasiconvexity (our Definition
\ref{def:O} is Hruska's (QC-5), and our Definition \ref{def:H} (QC-H)
is Hruska's (QC-3)), we obtain the
following useful fact.
\begin{corollary}\label{c:qcequiv}
  Let $G$ be relatively hyperbolic, relative to $\mc{P}$, and let $H$
  be a finitely generated subgroup of $G$.  The following are equivalent:
  \begin{enumerate}
  \item $H$ is a relatively quasiconvex subgroup of $G$, in the sense
    of Definition \ref{def:O}.
  \item $H$ is a relatively quasiconvex subgroup of $G$, in the sense
    of Definition \ref{def:AGM}.
  \end{enumerate}
\end{corollary}

\section{On extending the main result in the presence of torsion}\label{app:torsion}

In this section, we give some idea of the changes necessary to prove
Theorem \ref{thm:1} (and therefore Theorem \ref{t:c1}) in the presence
of torsion.  In this section, $G$ is a relatively hyperbolic group,
hyperbolic relative to a finite collection $\mathcal{P}$ of LERF and
slender subgroups, and $H$ is some relatively quasiconvex subgroup of
$G$. 

The main difference is that we must deal with the possibility that our
relatively quasiconvex subgroup has finite but non-trivial maximal parabolic
subgroups.  Since a finite subgroup of a relatively hyperbolic group,
may intersect arbitrary collections of parabolic subgroups, we have to
ignore these intersections.  This is already handled in the arguments
of Section \ref{sec:comb} by only amalgamating with parabolic
subgroups which have infinite intersection with $H$ to obtain the
fully quasiconvex subgroup $Q$.

In Section \ref{sec:fill}, it is necessary to modify the definition of
$H$--filling as follows:
\begin{definition}
  (Alternate definition in the presence of torsion.)  Let $H$ be a
  relatively quasiconvex subgroup of $G$.  A filling $G
\longrightarrow G(N_1, \dots , N_m)$ is an \emph{$H$--filling} if
whenever $H \cap P_i^g$ is infinite, $N_i^g \subset P_i^g \cap H$.
\end{definition}

With the new definition, we must check that the results from
\cite[Section 4.2]{AGM08} still hold.  (We do not know how to prove the
result about height from Section 4.3 of \cite{AGM08} in this more
general setting, but we do not need it for our argument.)  Examining
the proofs from \cite{AGM08}, the reader may check that
it suffices to extend the technical \cite[Lemma 4.2]{AGM08}. 

We sketch how to do
so briefly, for the experts:  In \cite[Lemma 4.2]{AGM08}, the hypothesis of an
$H$--filling is used to deduce the existence of a nontrivial element
of $H$ which is also in a conjugate of a filling kernel fixing a
certain horoball from the fact that a geodesic between
elements of $H$ penetrates that horoball deeply.  
The heart of the argument is showing that
if the geodesic penetrates the horoball deeply,
the intersection of $H$ with the horoball stabilizer is infinite.  In
the torsion-free setting, it suffices to show that the intersection is
nontrivial.  The proof in the presence of torsion is given in the
previous appendix as~\ref{l:shallowtorsion}.  With this proposition, 
one can prove the extended version of \cite[Lemma 4.2]{AGM08} in a straightforward manner,
choosing slightly different constants to take the constant $R$ from
Proposition~\ref{l:shallowtorsion} into account.

The proofs of Propositions 4.3 and 4.5 of \cite{AGM08} go through in exactly the same way, and we obtain
the same statement as Theorem \ref{thm:3} above, but with the new meaning of $H$--filling.  Using Osin's Dehn filling result in place of
Theorem \ref{thm:4}, the rest of the proof of Theorem \ref{thm:1} goes through as written, with the exception that each mention of a
condition of the form ``$A\cap B\neq\{1\}$'' for $A$ and $B$ subgroups of $G$ should be replaced by ``$A\cap B$ is infinite''.

\bibliographystyle{plain} 
\bibliography{Xbib}

\end{document}

%% file: fig02.pdf_t
\begin{picture}(0,0)%
\includegraphics{fig02.pdf}%
\end{picture}%
\setlength{\unitlength}{4144sp}%
\begingroup\makeatletter\ifx\SetFigFont\undefined%
\gdef\SetFigFont#1#2#3#4#5{%
  \reset@font\fontsize{#1}{#2pt}%
  \fontfamily{#3}\fontseries{#4}\fontshape{#5}%
  \selectfont}%
\fi\endgroup%
\begin{picture}(4347,1514)(3676,-4648)
\put(3861,-3979){\makebox(0,0)[lb]{\smash{{\SetFigFont{11}{13.2}{\rmdefault}{\mddefault}{\updefault}{\color[rgb]{0,0,0}$o$}%
}}}}
\put(3849,-4443){\makebox(0,0)[lb]{\smash{{\SetFigFont{11}{13.2}{\rmdefault}{\mddefault}{\updefault}{\color[rgb]{0,0,0}$u_{i-1}$}%
}}}}
\put(5889,-4443){\makebox(0,0)[lb]{\smash{{\SetFigFont{11}{13.2}{\rmdefault}{\mddefault}{\updefault}{\color[rgb]{0,0,0}$u_i$}%
}}}}
\put(6910,-3893){\makebox(0,0)[lb]{\smash{{\SetFigFont{11}{13.2}{\rmdefault}{\mddefault}{\updefault}{\color[rgb]{0,0,0}$v_i$}%
}}}}
\put(7905,-4416){\makebox(0,0)[lb]{\smash{{\SetFigFont{11}{13.2}{\rmdefault}{\mddefault}{\updefault}{\color[rgb]{0,0,0}$u_{i+1}$}%
}}}}
\put(4790,-3541){\makebox(0,0)[lb]{\smash{{\SetFigFont{11}{13.2}{\rmdefault}{\mddefault}{\updefault}{\color[rgb]{0,0,0}$a_iP_i$}%
}}}}
\put(6831,-3541){\makebox(0,0)[lb]{\smash{{\SetFigFont{11}{13.2}{\rmdefault}{\mddefault}{\updefault}{\color[rgb]{0,0,0}$a_{i+1}P_i$}%
}}}}
\put(3691,-3305){\makebox(0,0)[lb]{\smash{{\SetFigFont{11}{13.2}{\rmdefault}{\mddefault}{\updefault}{\color[rgb]{0,0,0}$\Gamma(G, \mathcal{P}, S)$}%
}}}}
\put(4861,-3886){\makebox(0,0)[lb]{\smash{{\SetFigFont{11}{13.2}{\rmdefault}{\mddefault}{\updefault}{\color[rgb]{0,0,0}$v_{i-1}$}%
}}}}
\end{picture}%

%% file: fig01.pdf_t
\begin{picture}(0,0)%
\includegraphics{fig01.pdf}%
\end{picture}%
\setlength{\unitlength}{4144sp}%
\begingroup\makeatletter\ifx\SetFigFont\undefined%
\gdef\SetFigFont#1#2#3#4#5{%
  \reset@font\fontsize{#1}{#2pt}%
  \fontfamily{#3}\fontseries{#4}\fontshape{#5}%
  \selectfont}%
\fi\endgroup%
\begin{picture}(4436,3562)(3811,-4937)
\put(4534,-4468){\makebox(0,0)[lb]{\smash{{\SetFigFont{12}{14.4}{\rmdefault}{\mddefault}{\updefault}{\color[rgb]{0,0,0}$1$}%
}}}}
\put(3984,-3445){\makebox(0,0)[lb]{\smash{{\SetFigFont{12}{14.4}{\rmdefault}{\mddefault}{\updefault}{\color[rgb]{0,0,0}$f$}%
}}}}
\put(7131,-1636){\makebox(0,0)[lb]{\smash{{\SetFigFont{12}{14.4}{\rmdefault}{\mddefault}{\updefault}{\color[rgb]{0,0,0}$fg$}%
}}}}
\put(7760,-2659){\makebox(0,0)[lb]{\smash{{\SetFigFont{12}{14.4}{\rmdefault}{\mddefault}{\updefault}{\color[rgb]{0,0,0}$fgf^{-1}$}%
}}}}
\put(8232,-4390){\makebox(0,0)[lb]{\smash{{\SetFigFont{12}{14.4}{\rmdefault}{\mddefault}{\updefault}{\color[rgb]{0,0,0}$g$}%
}}}}
\put(3826,-1558){\makebox(0,0)[lb]{\smash{{\SetFigFont{12}{14.4}{\rmdefault}{\mddefault}{\updefault}{\color[rgb]{0,0,0}$\Gamma(G, \mathcal{P}, S)$ }%
}}}}
\put(4771,-3052){\makebox(0,0)[lb]{\smash{{\SetFigFont{12}{14.4}{\rmdefault}{\mddefault}{\updefault}{\color[rgb]{0,0,0}$r$}%
}}}}
\put(5242,-3760){\makebox(0,0)[lb]{\smash{{\SetFigFont{12}{14.4}{\rmdefault}{\mddefault}{\updefault}{\color[rgb]{0,0,0}$q$}%
}}}}
\put(5400,-4547){\makebox(0,0)[lb]{\smash{{\SetFigFont{12}{14.4}{\rmdefault}{\mddefault}{\updefault}{\color[rgb]{0,0,0}$p$}%
}}}}
\put(6344,-4547){\makebox(0,0)[lb]{\smash{{\SetFigFont{12}{14.4}{\rmdefault}{\mddefault}{\updefault}{\color[rgb]{0,0,0}$w$ }%
}}}}
\put(5804,-3061){\makebox(0,0)[lb]{\smash{{\SetFigFont{12}{14.4}{\rmdefault}{\mddefault}{\updefault}{\color[rgb]{0,0,0}$u$ }%
}}}}
\put(6712,-3916){\makebox(0,0)[lb]{\smash{{\SetFigFont{12}{14.4}{\rmdefault}{\mddefault}{\updefault}{\color[rgb]{0,0,0}$bP_i$ }%
}}}}
\put(5637,-2667){\makebox(0,0)[lb]{\smash{{\SetFigFont{12}{14.4}{\rmdefault}{\mddefault}{\updefault}{\color[rgb]{0,0,0}$v$ }%
}}}}
\put(5047,-2364){\makebox(0,0)[lb]{\smash{{\SetFigFont{12}{14.4}{\rmdefault}{\mddefault}{\updefault}{\color[rgb]{0,0,0}$aP_i$ }%
}}}}
\end{picture}%